\documentclass[runningheads]{llncs}
\usepackage{graphicx}
\usepackage{wrapfig, lipsum}
\usepackage{tikz}
\usepackage{amsmath}
\usepackage{amsfonts}
\usepackage{float}
\usepackage[section]{placeins}
\usepackage{cleveref}

\begin{document}
\title{Reduced Basis Methods for Efficient Simulation of a Rigid Robot Hand Interacting with Soft Tissue\thanks{The research leading to this publication has received funding from the German Research Foundation (DFG) as part of the International Research Training Group ``Soft Tissue Robotics'' (GRK 2198/1) and under Germany's Excellence Strategy - EXC 2075 – 390740016. We acknowledge the support by the Stuttgart Center for Simulation Science (SimTech).}}  
\titlerunning{Efficient Reduced Basis for Soft Tissue Robotics}
% If the paper title is too long for the running head, you can set an abbreviated paper title here
\author{Shahnewaz Shuva\inst{1} \and Patrick Buchfink\inst{1} \and
Oliver R\"ohrle\inst{2} \and
Bernard Haasdonk\inst{1}}
%\orcidID{0000-1111-2222-3333}
\authorrunning{S. Shuva et al.}
% First names are abbreviated in the running head.
% If there are more than two authors, 'et al.' is used.
\institute{Institute of Applied Analysis and Numerical Simulation\\ University of Stuttgart, Germany \\
\email{\{shuvasz, buchfipk, haasdonk\}@mathematik.uni-stuttgart.de}\\
 \and
Institute for Modeling and Simulation of Biomechanical Systems\\ University of Stuttgart, Germany \\
\email{roehrle@simtech.uni-stuttgart.de}}
\maketitle              % typeset the header of the contribution
\begin{abstract}
We present efficient reduced basis (RB) methods for the simulation of a coupled problem consisting of a rigid robot hand interacting with soft tissue material. The soft tissue is modeled by the linear elasticity equation and discretized with the Finite Element Method. We look at two different scenarios: (i) the forward simulation and (ii) a feedback control formulation of the model. In both cases, large-scale systems of equations appear, which need to be solved in real-time. This is essential in practice for the implementation in a real robot. For the feedback-scenario, we encounter a high-dimensional Algebraic Riccati Equation (ARE) in the context of the linear quadratic regulator. To overcome the real-time constraint by significantly reducing the computational complexity, we use several structure-preserving and non-structure-preserving reduction methods. These include reduced basis techniques based on the Proper Orthogonal Decomposition. For the ARE, we compute a low-rank-factor and hence solve a low-dimensional ARE instead of solving a full dimensional problem. Numerical examples for both cases (i) and (ii) are provided. These illustrate the approximation quality of the reduced solution and speedup factors of the different reduction approaches.
\end{abstract}
\section{Introduction}
\begin{figure}[!htb]
\centering
\begin{tikzpicture}[baseline={(current bounding box.center)},scale=.4, grid/.style={very thin,gray} ]         
%\draw[thick,red] (-.5, 5) .. controls (1.5,7)and(4,4) .. (10, 5);
%\node[] at (3,6) {Obstacle for $\Sigma_{s}$};
\draw [thick](0,1) -- (-1,1) -- (-1, 4) -- (3,4) -- (3,1) -- (2,1) -- (2, 2.5) --(0, 2.5) -- (0,1); % robot handle
\node[] at (0, 3.5) {$\Sigma_{s}$};
\draw [blue](0,1) -- (0,2) -- (2,2) -- (2,1);
\draw [gray](0,1) -- (2,1);
\draw [gray](0,1.5) -- (2,1.5);
\draw [gray](.5,2) -- (.5,1);
\draw [gray](1,2) -- (1,1);
\draw [gray](1.5,1) -- (1.5,2);
\draw [gray](0.5,1) .. controls (0.6,.5) .. (.5,-.025);
\draw [gray](1,1) .. controls (1.1,.5) .. (1,-.05);
\draw [gray](1.5,1) .. controls (1.6,.5) .. (1.5,-.025);
\draw [green](0.05,.5) .. controls (1,.4) .. (2.1,0.49);
\draw [red, thick](0,0) .. controls (1,-.1) .. (2,0);
\draw [red,thick](0,1) .. controls (0.1,.5) .. (0,0);
\draw [red, thick](0,0) .. controls (1,-.1) .. (2,0);
\draw [red, thick] (2,0) .. controls (2.1,.5) .. (2,1);
\node[] at (1,1) {$\Sigma_{e}$};
\node[] at (1,4.5) {$\Sigma$};
\filldraw[black] (1,3.5) circle (2pt) node[anchor=west]{};
\draw (1,3.5) .. controls (3,3) and (4,5) .. (9,4); %pathline
\filldraw[black] (9, 4) circle (2pt) node[anchor=west]{};
\node[] at (4,3.5) {$u(t)$};
\draw [gray](8,0.5) -- (10,0.5)  -- (10,2.5)  -- (8,2.5) -- (8,0.5);
\node[scale =.5] at (9,2) {Target};
\node[scale = .5] at (9,1.5) {Position}; 
\draw [ -latex] (-.1,.25) -- (-.4,.25);
\draw [ -latex] (2.1,.25) -- (2.4,.25);
\end{tikzpicture}
\caption{Schematic view of gripper soft tissue system with target position.}
\label{fig:figure1}
\end{figure}
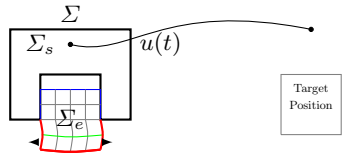

The ability to manipulate deformable objects using robots has diverse applications with enormous economical benefits. The applications include food industry, medical sectors, automobile industry, soft material processing industry and many more. Although grasping and manipulation of rigid objects by robots is a mature field in robotics, with over three decades of works, the study of deformable objects has not been as extensive in the robotics community \cite{stateofart_1}. Here, we model the elastic object using existing linear elastic theory with two characterization parameters, namely, first and second Lam\'e parameter. The two most challenging and frequently studied manipulation tasks on planar deformable objects are grasping and controlling its deformations \cite{grasping}. Grasping an object consists of positioning the end effectors of a robot hand on the object to lift and hold it in the air, which involves the challenges of slipping off and non-linear contact mechanics.
For the sake of simplicity, we assume the soft tissue material to be attached to the rigid robot's end effectors. We are interested in controlling the object using a feedback controller and transfer it to a target position. 
Here, we focus on two different scenarios, (i) the forward simulation of the coupled problem, where the robot hand along with the soft tissue material follows a prescribed trajectory, and (ii) a feedback control such that the robot hand along with the soft tissue material cost-optimally reaches a target position and then stabilizes. In Fig.~\ref{fig:figure1}, a schematic view is provided.
%====================================================================
After the discretization of both problems, large-scale systems of equations appear. Another challenging issue besides modeling and controlling is that the simulation of the model has to be finished in real-time because in a modern robotics software, like Franka Control Interface \cite{libfranka}, all real-time loops have to be finished within 1 ms. %Each controller needs to return values within $1$ kHz. 
A feasible solution of this is to simulate a reduced model instead of a full order model. We developed and applied different structure-preserving and non-structure-preserving methods to maintain specific structures, e.g.\ block structure in the reduced system.
The paper is organized in the following way. In Section $2$, we discuss the modeling based on a spatially continuous formulation, which is followed by a large-scale spatially discrete problem in the context of the forward and LQR control problem. In Section $3$, we discuss different reduced basis methods for both cases. Next in Section $4$, we apply the MOR techniques in two numerical examples. Finally, we compare the approximation quality of the reduced solutions and provide a comparison of the execution time in contrast to the full order model.
\vspace*{-.2cm}  
% =============================================================
\section{Model Problems}
We model the soft tissue material by the following time-dependent linear elasticity partial differential equation (PDE). We assume a two-dimensional spatial domain $\Omega_{e}:= (0, 1)\times (0,2) \subset \mathbb{R}^2$ with boundary, $\Gamma = \Gamma_D \cup \Gamma_N$ where $\Gamma_{D} = \{0,1\} \times [1.5, 2]$ denotes the Dirichlet boundary and $\Gamma_{N} = \Gamma \setminus \Gamma_D$ denotes the Neumann boundary.
The governing equation for the displacement field is
\begin{align*}
\nabla	\cdot \boldsymbol{\sigma} + \boldsymbol{F} &= \rho \ddot{\boldsymbol{q}}^{c} \qquad \text{in} \qquad \Omega_e\times (0,T)
\end{align*}
with suitable initial and boundary conditions where $\boldsymbol{q}^{c}$ is the displacement vector and $\boldsymbol{\sigma}$ the Cauchy stress tensor. $\boldsymbol{F}$ is the external body force %(gravity, centrifugal force or similar) 
per unit volume and $\rho$ is the mass density. The constitutive equation is $\boldsymbol{\sigma} = C:\boldsymbol{\epsilon}$ according to Hooke's law for elastic material, where $C$ is a fourth-order tensor and $:$ is the contraction operator. For isotropic and homogeneous media holds $ \boldsymbol{\sigma} = 2 \mu \boldsymbol{\epsilon} +\lambda \textit{tr}(\boldsymbol{\epsilon}) \boldsymbol{I} $ where $\lambda$, $\mu$ are the Lam\'e parameters.
The strain is expressed in terms of the gradient of the displacement with $\boldsymbol{\epsilon} = \frac{1}{2}(\nabla \boldsymbol{q}^{c} + {(\nabla \boldsymbol{q}^{c})}^T)$.
\subsection{Forward Model}
After discretization %of the linear elasticity PDE 
using the Finite Element Method (FEM), % and coupling with a solid gripper model
we get the following second order system of differential equations,
\begin{align}
 M(\rho) \ddot{q} (t) + K (\mathbf{\mu, \lambda}) q (t) = f(t) + B_{u} u(t) \label{eq:1}
\end{align}
with the parameters $\rho, \lambda, \mu \in \mathbb{R}$ and the state vector $q \in \mathbb{R}^{n}$ which is decomposed into the displacement vector $q_s \in \mathbb{R}^{n_s}$ of the solid robot hand and the displacement vector $q_e \in \mathbb{R}^{n_e}$ of the soft tissue material, i.e.\ $q := [q_s, q_e]^{T} \in \mathbb{R}^{n}$, $n := n_s + n_e$. The mass and stiffness matrix $M, K \in \mathbb{R}^{n \times n}$, the influence matrix $B_u \in \mathbb{R}^{n \times m}$, the body force $f(t) =  \begin{bmatrix} f_{s}^T(t) & f_{e}^T(t) \end{bmatrix}^T \in \mathbb{R}^n$ and input vector $u(t) = \begin{bmatrix} u_{s}^T(t) & u_{e}^T(t) \end{bmatrix}^T \in \mathbb{R}^n$ in the above equations are as follows:
$ M = \begin{bmatrix}
M_{ss} & M_{se};
M_{es} & M_{ee}
\end{bmatrix} $, 
$K = \begin{bmatrix}
K_{ss} & K_{se};
K_{es} & K_{ee}
\end{bmatrix} $, $B_u = \text{blkdiag}(B_{us},B_{ue})$. 
Here the block matrices with the subscript $se$ or $es$ refer to the coupled matrix coefficient between the vectors $q_s$ and $q_e$ and the $;$ indicates a new block row. With all non-zero block matrices, the problem can be interpreted as a two-way coupled problem \cite{blockstrucpreserv_1}. For simplification, we assume that the mass of the solid system is large enough compared to the mass of the elastic system, such that the influence of the motion of the elastic system on the solid body is negligible but vice versa is relevant, i.e.\ $M_{se}$ is chosen as a zero matrix and $M_{es}$ is a non-zero matrix. Since we consider the robot hand as a rigid body, $K_{ss}$ is set to zero. The contact zone follows the motion of the solid body, so we ignore $K_{se}$. The values of the displacement vectors on the Dirichlet boundary of the elastic body are determined by the displacement vectors of the solid body $q_s$, i.e.\ $K_{es}$ is a non-zero matrix. With these assumptions, we obtain the following one-way coupled model
\begin{align}
\begin{bmatrix}
M_{ss} & 0 \\
M_{es} & M_{ee}
\end{bmatrix} \begin{bmatrix}
\ddot{q}_{s}(t) \\ \ddot{q}_{e}(t) 
\end{bmatrix} + \begin{bmatrix}
0 & 0 \\
K_{es} & K_{ee}
\end{bmatrix} \begin{bmatrix}
q_{s}(t) \\ q_{e}(t)  
\end{bmatrix} = \begin{bmatrix}
f_s(t) \\  f_e(t) \end{bmatrix} +
\begin{bmatrix} B_{us} & \\ & 0 \end{bmatrix} \begin{bmatrix} u_{s}(t) \\  u_{e}(t) \end{bmatrix}_. \label{Eq:12} 
\end{align}
Here, the elasticity part of the state vector is coupled with the acceleration and displacement of the solid state vector. For solving Eq.~\ref{Eq:12}, we transformed it into system of first order differential equations, where the state vector is composed of the displacement $q(t)$ and the velocity $v(t)=\dot{q}(t)$.
\subsection{Linear Quadratic Regulator (LQR)}
For feedback control, we formulate a LQR problem and assume the weighing matrices $Q\in \mathbb{R}^{p \times p}, \; p \leq 2n $ to be symmetric positive semi-definite and $R\in \mathbb{R}^{m \times m}$ to be symmetric positive definite. We define the quadratic cost functional
\begin{align*}
%J(u) & = \int_{0}^{\infty} [y(t)^{T}Q \, y(t) + u(t)^{T}R\,u(t)] dt \nonumber \\
J(u) = \int_{0}^{\infty} [x(t)^{T}C^TQ^{1/2} Q^{1/2} \,C x(t) + u(t)^{T}R\,u(t)] dt \nonumber
\end{align*}
with $ E \dot{x}(t) = A x(t) + B u(t),\: y(t) = Cx(t), \: x(0) = x_0 $, time $t \in \mathbb{R}_{+} $, the state $x(t) = \begin{bmatrix}
q_{s}^T(t) & \dot{q}^T_{s}(t) & q_{e}^T(t) & \dot{q}_{e}^T(t)
\end{bmatrix}^T \in \mathbb{R}^{2n}$, the input $u(t) \in \mathbb{R}^m$, the output $y(t) \in \mathbb{R}^{p}$ and system matrices $E,A \in \mathbb{R}^{2n \times 2n}, B \in \mathbb{R}^{2n \times m}$, $C \in \mathbb{R}^{p \times 2n}$ with
$E=  \begin{bmatrix} E_{ss} & 0 \\ E_{es} & E_{ee} \end{bmatrix},\;
A= \begin{bmatrix} A_{ss} & 0\\A_{es} & A_{ee}\end{bmatrix},\;
B  = \begin{bmatrix} 0 \\B_{us} \\ 0 \end{bmatrix},
\setcounter{MaxMatrixCols}{20}
 C = \begin{bmatrix}
 0 &0 & 0 & 0 & 1 & 0 &\hdots& 0 & 0 & 0 & \hdots & 0\\
 0 & 0 & 0 & 0 & 0 & 0 &\hdots& 0 & 1 & 0 & \hdots & 0
\end{bmatrix}.$\\
If we assume $(E, A, B)$ to be stabilizable and $\begin{pmatrix} E, A, Q^{1/2}C\end{pmatrix}$ to be detectable, the optimal control problem (OCP) possesses a unique solution $u(t) = K_f x(t)$ with the feedback gain matrix $K_f = -R^{-1}B^TPE \in \mathbb{R}^{m \times 2n}$. Here, the symmetric positive semi-definite matrix $P \in \mathbb{R}^{2n \times 2n}$ is the unique stabilizing solution of the Generalized Algebraic Riccati Equation (ARE) \cite{lqr_are_1}
\begin{align}
E^{T}PA + A^{T}PE  - E^{T}PBR^{-1}B^{T}PE + C^{T}QC = 0 \label{Eq:2}
\end{align}
Note that this is a large-scale problem with $4n^2$ unknowns and therefore of quadratic complexity.
\vspace*{-.3cm}
\section{Reduced Basis Methods}
%\vspace*{-.3cm}
The RB-method \cite{rbmatlab_bookchapt} approximates the solution in a low-dimensional subspace that is constructed from solutions of the large-scale problems (Eq.\ \ref{eq:1} and \ref{Eq:2}). To construct a reduced basis matrix $V := \mathrm{span}\{ \psi_1, \hdots \psi_{N_V} \} \in \mathbb{R}^{2n \times N_V}$ with $N_V$ basis functions, we define a snapshot matrix $ X := [x_1, \hdots, x_{n_{sn}}] \in \mathbb{R}^{2n \times n_{sn}}$ with $1\leq i \leq n_{sn}$ so-called snapshots $x_i$, which are acquired by solving the corresponding large-scale problem. For different reduction techniques, we use different splittings of the snapshot matrix: a split $X_e \in \mathbb{R}^{2n_e \times n_{sn}}$, $X_s \in \mathbb{R}^{2n_s \times n_{sn}}$ in elastic and solid part, a split $X_q, X_v \in \mathbb{R}^{n \times n_{sn}}$ in velocity and displacement and a combination of both splits $X_{s_q}, X_{s_v} \in \mathbb{R}^{n_s \times n_{sn}}, X_{e_q}, X_{e_v} \in \mathbb{R}^{n_e \times n_{sn}}$. The reduction techniques are based on the so-called Proper Orthogonal Decomposition (POD). It chooses the reduced basis as the first $k$ left-singular vectors of the snapshot matrix $X$ which we denote with $\text{POD}_k(X)$. In practice, the POD basis is computed via the truncated Singular Value Decomposition of $X$. In the following, the methods are classified as non-structure-preserving and structure-preserving techniques. 
%\vspace*{-.4cm}
\subsection{MOR for forward problem}
For the model reduction, we apply Galerkin projection with a suitable basis matrix $V$. The non-structure-preserving methods are the following:
\let\labelitemi\labelitemii
\begin{itemize}
\item \textbf{Global POD}: Here, we implement a global POD, i.e.\ perform a singular value decomposition (SVD) of the snapshot matrix, $V = \mathrm{POD}_k(X)$.
\item \textbf{Componentwise POD}: In this technique, we apply the POD algorithm for displacement and velocity snapshots separately, i.e.\ $V_1 = \mathrm{POD}_{k_1}(X_q)$, $V_2 = \mathrm{POD}_{k_2}(X_v)$ and then construct $V = \text{blkdiag}(V_1, V_2)$.
\end{itemize}
As structure-preserving techniques, we apply the following methods:
\begin{itemize}
\item \textbf{POD with fixed solid modes}: The reduced basis space is constructed in the following way: We apply POD only in the elastic part of the snapshots, i.e.\ $V_e = \mathrm{POD}_{k_e}(X_e)$ and hence construct $V = \mathrm{blkdiag} (I_{2n_s}, V_e)$ where $I_{2n_s} \in \mathbb{R}^{2n_s \times 2n_s}$ is the identity matrix. The unreduced solid modes are still computationally efficient since $n_s \ll n_e$.
\item \textbf{Componentwise POD with fixed solid modes}: Here we construct $V= \text{blkdiag}(I_{2n_s}, V_{e_q}, V_{e_v})$, with $V_{e_q} = \mathrm{POD}_{k_{e_q}}(X_{e_q})$, $V_{e_v} = \mathrm{POD}_{k_{e_v}}(X_{e_v})$.
\item \textbf{Global Proper Symplectic Decomposition(GPSD)}: As our system for $u=0$ is a Hamiltonian system, we construct a so-called orthosymplectic basis matrix. As basis generation technique, we use the so-called PSD Complex SVD based on the SVD of a modified snapshot matrix $Y = [X\; JX]$, where $J$ is the corresponding Poisson matrix. For more details, we refer to \cite{symplectic_1}.
%For details we refer to Prop. 1 in  We augment the snapshot matrix with ``rotated" snapshots to $Y = \begin{bmatrix}X  &  \mathbb{J}X \end{bmatrix}$ with Poisson matrix $\mathbb{J}$ which is a suitable skew-symmetric block matrix. We assume the $2k$ such that there is a gap in the singular values of $Y$, i.e., $\sigma_{2k}(Y) > \sigma_{2k+1}(Y)$. Then, GPSD in the set of symplectic, orthogonal reduced order basis is equivalent to solve minimization problem
%$\underset{V \in \mathbb{R}^{2n,2k}}{\text{min}} \|(I_{2n} - VV^T) Y \|^2_F \quad \text{s.t.}\, V^TV = I_{2k}.$ 
\item \textbf{PSD with fixed solid modes}: Here, we apply the PSD only on the elastic part.
\end{itemize}
\vspace*{-.3cm}
\subsection{RB-ARE for LQR}
In practice, the dimension of the state space and hence of $P$ is typically very high. The dimension $m$ of the input $u(t)$ and $p$ of the output $y(t)$ is much smaller than the number of states $2n$. The computation of $P$ is often very expensive or even impossible \cite{lqr_are_2}. Instead of determining $P$, we compute $ Z Z^T \approx P $ with low-rank-factor $Z \in \mathbb{R}^{2n \times N}$ of low-rank $N \ll 2n$. The methods without fixing the solid modes include:
\begin{itemize}
\item \textbf{POD of ARE solution $P$}: We construct the reduced basis matrix using the POD of the solution of ARE, i.e.\ $V = \mathrm{POD}_k(P)$. 
\item \textbf{Weighted POD}: We assume a weighting matrix $W = \frac{1}{2}(E + E^T)$ and apply a weighted POD \cite{pod}, i.e.\ $V := W^{-\frac{1}{2}}\mathrm{POD}_k(W^{\frac{1}{2}}P)$.
\end{itemize}
The method with fixed solid modes are the following:
\begin{itemize}
\item \textbf{POD with fixed solid modes}: We decompose the solution $P$ of the ARE into $P_{ss}$, $P_{se}$, $P_{es}$ and $P_{ee}$ and apply $ V_e := \mathrm{POD}_k 
\begin{pmatrix}
\begin{bmatrix}
P_{es}, P_{ee}
\end{bmatrix}
\end{pmatrix}$ and construct $V := \mathrm{blkdiag}(I_{2n_s}, V_e).$
%\begin{bmatrix}
%P_{ss} & P_{se}\\ P_{es} & P_{ee}
%\end{bmatrix}$
\item \textbf{Componentwise POD}: We decompose blocks of $P$ in displacement and velocity and use $\mathrm{POD}_{k_1}\begin{pmatrix}
\begin{bmatrix}
P_{es_{11}}, P_{ee_{11}},P_{ee_{12}}
\end{bmatrix}
\end{pmatrix}$ and $\mathrm{POD}_{k_2}\begin{pmatrix}
\begin{bmatrix}
P_{es_{21}}, P_{ee_{21}},P_{ee_{22}}
\end{bmatrix}
\end{pmatrix}$ for displacement and velocity part separately to construct V.
\item \textbf{POD of decomposed $P$}: We use $V_1 := \mathrm{POD}_{k_1}(  
\begin{bmatrix}
P_{ee_{11}}, P_{ee_{12}}, P_{ee_{21}}, P_{ee_{22}} 
\end{bmatrix})$ where $P_{ee_{11}}, P_{ee_{12}}, P_{ee_{21}}, P_{ee_{22}} \in \mathbb{R}^{n_e \times n_e}$ with $V:= \mathrm{blkdiag}(I_{2n_s}, V_1, V_1).$
\end{itemize} 
We introduce the reduced basis approximation $\hat{P} = VP_{N}V^{T}$, where $N \ll 2n $ and $P_N \in \mathbb{R}^{N \times N}$ is the solution of the reduced ARE
\begin{align*}
E_N ^{T}P_N A_N  + A_N ^{T}P_N E_N   - E_N ^{T}P_N B_N R ^{-1}B_N ^{T}P_N E_N  + C_N ^{T}QC_N  = 0. 
\end{align*}
$V$ is constructed by one of the above mentioned methods (one parameter setting), $A_N = V^T A V$, $E_N = V^T E V$, $B_N = V^T B$ and $C_N = CV$.
\vspace*{-.3cm}
\section{Numerical Examples}
We present our numerical experiments for both, the forward and the LQR problem. As a time integrator, we use the implicit mid-point rule with $n_t = 600$ time steps. We choose ranges of the first and second Lam\'e parameter, $\lambda \in [30, 500]$ and $\mu \in [20, 500]$. We consider one trajectory for a particular value of the parameters ($\lambda =50$, $\mu = 50$) which results in $n_{sn} = n_{t}$ snapshots. The initial state is $x(0) = 0 \in \mathbb{R}^{2n}$. As a target position we consider $\bar{x}(T) = \begin{bmatrix} 5 &5 & \hdots & 3 & 3 & \hdots \end{bmatrix}^T$ with the final time $T = 3$ and $300$ for the forward and the LQR problem, respectively. For the simulation and reduction of the model, we use the software package RBmatlab \cite{rbmatlab_bookchapt}.
\subsection{Forward Problem}
\begin{figure}[H]
\centering
\includegraphics[width=.33\textwidth]{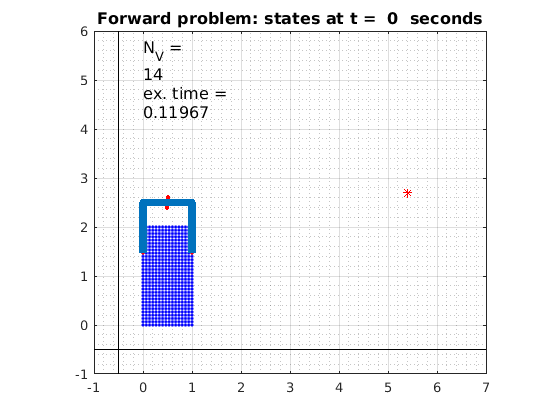}\hfill
\includegraphics[width=.33\textwidth]{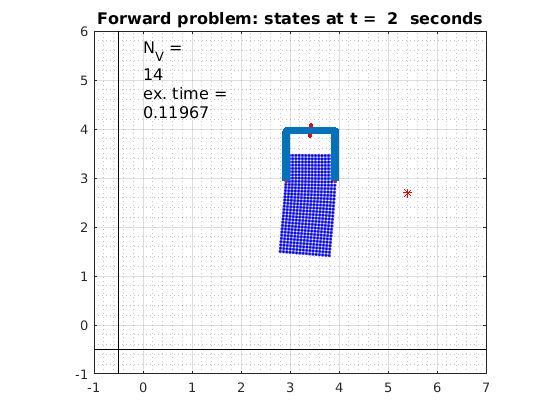}\hfill
\includegraphics[width=.33\textwidth]{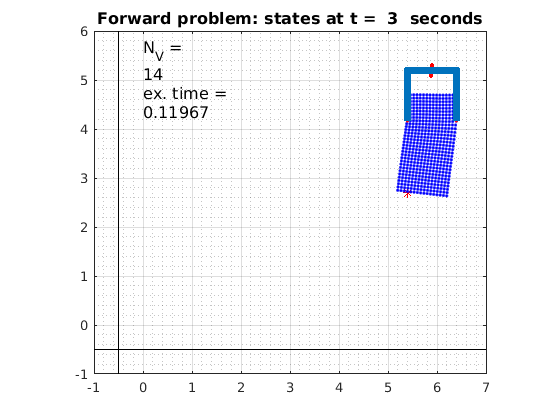}
\caption{Forward simulation of the reduced coupled problem}
\label{Fig:fig2}
\end{figure}

\begin{figure}[!htp]
\centering
  \includegraphics[width=.5\textwidth]{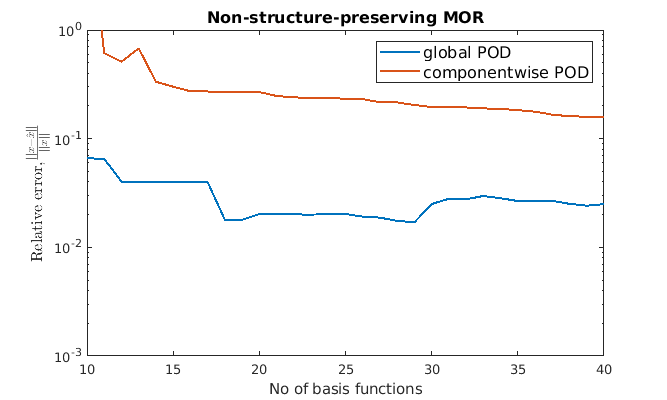}\hfill
  \includegraphics[width=.5\textwidth]{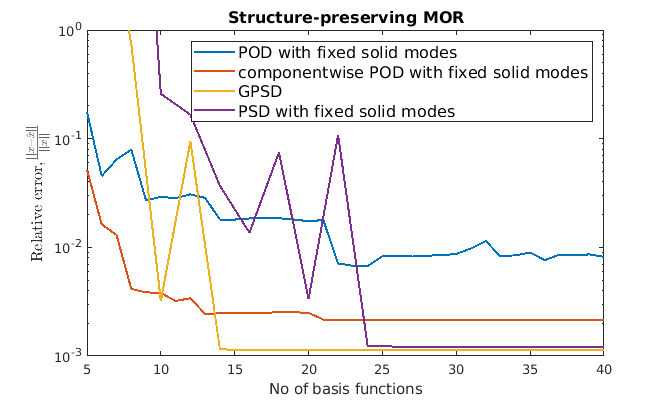}\hfill
  \caption{Approximation quality of reduced solutions.}
  \label{fig:fig3}
\end{figure}

We determine $N_V$ to include $99.9\%$ energy of the POD functional. An example of a reduced forward problem is illustrated in Fig.\ \ref{Fig:fig2} using the PSD method. To measure the quality of the reduced solution we compute the relative error $\frac{ \| X-\hat{X} \|_F}{\| X\|_F}$ (see Fig.\ \ref{fig:fig3}), where $X \in \mathbb{R}^{2n \times n_{sn}}$ is the snapshot matrix introduced above, $ \hat{X}: = VX_r$ is the reconstructed solution, $X_r \in \mathbb{R}^{N_V \times n_{sn}}$ is the matrix gathering the reduced solution time-instances as columns and $\|.\|_F$ is the Frobenius norm.
It shows that the structure-preserving methods are better than the non-structure-preserving methods. 
The best technique is the PSD without fixed solid modes which results in an error of $\approx 10^{-3}$ for the reduced system when using at least $N_V=14$ basis functions.
The execution time for the different methods is between $0.097s$ and $0.101s$ depending on the number of basis functions. The PSD with $N_V=14$ requires $0.097s$ which is much less compared to the execution time $16.22s$ for the high-dimensional solution ($2n=1916$). For a size of the reduced models larger than 25, we do not get any improvements in accuracy because the singular values do not show a sharp decay after that.
%\vspace*{-1cm}
\subsection{LQR Problem}
\begin{figure}[H]
\centering
\includegraphics[width=.33\textwidth]{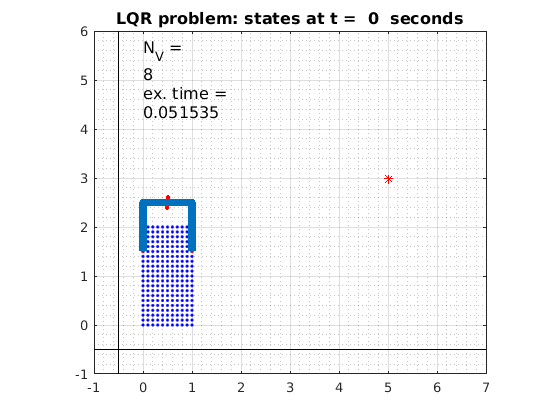}\hfill
\includegraphics[width=.33\textwidth]{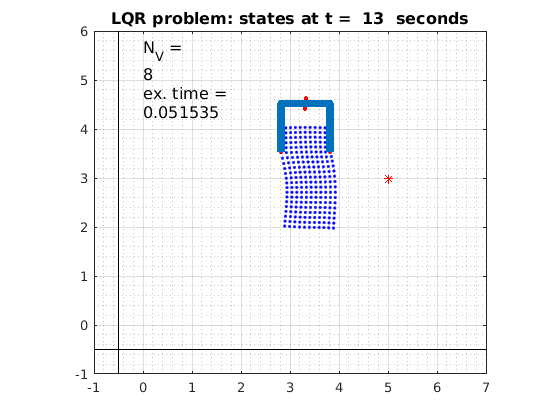}\hfill
\includegraphics[width=.33\textwidth]{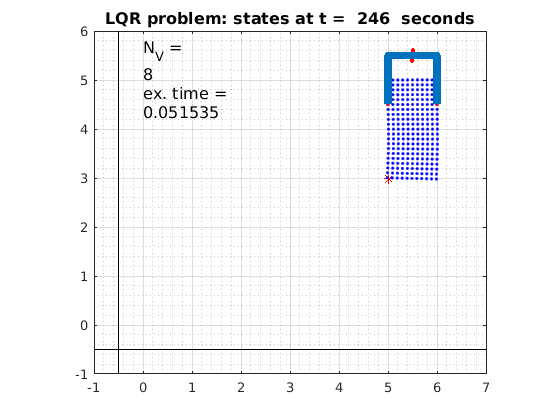}
\caption{Simulation of the reduced LQR problem}
\label{fig:fig4}
\end{figure}

\begin{figure}[htb]
\centering
  \includegraphics[width=.5\textwidth]{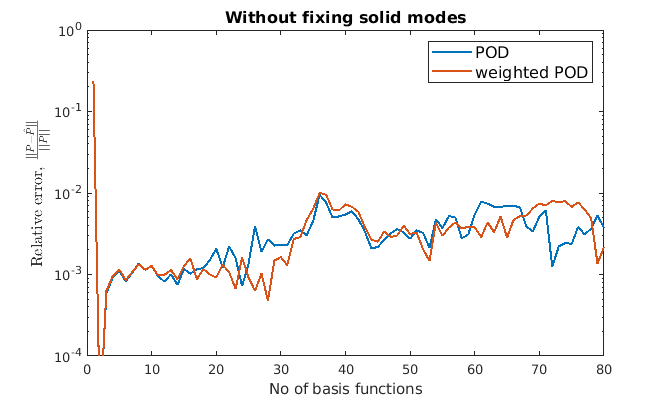}\hfill
  \includegraphics[width=.5\textwidth]{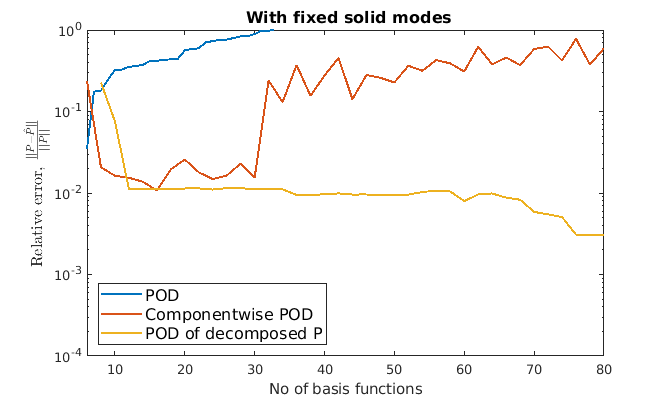}\hfill
  \vspace*{-.01cm}
  \caption{Relative error for the controller-scenario. Errors above $100\%$ are not depicted.}
  \label{fig:fig5}
\end{figure}

We first provide the simulation using the reduced controller in Fig.~\ref{fig:fig4} using $V$ from the global POD method. It shows that the reduced solution is able to compute a reconstructed stabilizing solution $\hat{P}$ which drives the soft tissue to its target position. 
%We applied the MOR methods which we listed in the Section $3.2$. Here we compute the quality of the reduced ARE solution. We first provide the simulation using the reduced controller. 
\vspace*{-.5cm}
We applied the MOR methods listed in Section $3.2$. The relative error of the reduced ARE solution is shown in Fig.~\ref{fig:fig5} which indicates that the methods without fixed solid modes provide better results for smaller numbers of basis functions. We also see that the computation time for Eq.~\ref{Eq:2} ($2n = 880$) requires $84.53s$ which we can reduce to a runtime between $0.051s$ and $0.057s$ depending on the number of basis functions, e.g.\ $0.053s$ for the global POD with $N_V = 8$.
\vspace*{-2cm}
\section{Conclusion and Outlook}
\vspace*{-.1cm}
We compared the quality and speedup of different basis generation techniques for a large-scale soft tissue model in the context of a forward problem and a feedback control problem. In the forward problem, we see that the PSD-based methods provide better results in terms of accuracy than the POD-based methods. For two example reduced models with a good accuracy, we are able to achieve speedup factors of 167 and 1600, respectively, which motivates us to apply these methods in parametric problems for multi-query scenarios. The future work includes adding an obstacle using state constraints in order to look at a more general class of control problems, greedy basis generation for the multi-query case and modeling using a non-linear elastic material law. 

\end{document}